\documentclass[11pt,reqno,letterpaper]{amsart}
\usepackage[T1]{fontenc}
\usepackage{lmodern}
\usepackage{amsmath,amssymb,amsthm}
\usepackage{microtype}
\usepackage[margin=28mm]{geometry}
\usepackage[hidelinks,hypertexnames=false]{hyperref}

\newtheorem{theorem}{Theorem}[section]
\newtheorem{corollary}[theorem]{Corollary}
\newtheorem{proposition}[theorem]{Proposition}
\newtheorem{lemma}[theorem]{Lemma}
\theoremstyle{remark}
\newtheorem{remark}[theorem]{Remark}
\newcommand{\R}{\mathbb R}
\newcommand{\MM}{\widetilde{\mathfrak M}}
\newcommand{\dd}{\,\mathrm d}
\newcommand{\PP}{\mathcal P_0}
\DeclareMathOperator{\Lip}{Lip}
\numberwithin{equation}{section}

\title[Nonlinear selection and vanishing discount]{Nonlinear selection in vanishing discount problems for Hamilton--Jacobi equations}
\author{Panrui Ni}
\address[P. Ni]{Department of Mathematics, Faculty of Science and Engineering, Waseda University, 3-4-1 Okubo, Shinjuku-ku, Tokyo 169-8555, Japan; Shanghai Center for Mathematical Sciences, Fudan University, Shanghai 200438, China}
\email{panruini@gmail.com}
\author{Maxime Zavidovique}
\address[M. Zavidovique]{Sorbonne Universit\'e, Universit\'e de Paris Cit\'e, CNRS, Institut de Math\'ematiques de Jussieu-Paris Rive Gauche, Paris 75005, France}
\email{mzavidovi@imj-prg.fr}
\date{}

\begin{document}

\begin{abstract}
We study the vanishing discount problem
$G(x,Du_\lambda,\lambda u_\lambda)=c_0$, where $G$ is convex and
superlinear in the momentum and strictly increasing in its last variable.
Assuming a leading-order expansion of the Lagrangian of order $q>0$,
we prove uniform convergence and characterize the selected solution by
nonlinear inequalities against Mather measures. In particular, our
results cover discount terms whose first derivative at the origin
vanishes, such as the cubic discount, beyond existing convergence
results based on first-order nondegeneracy. Our proof uses truncated
calibrated curves and time averages, without the exponentially
weighted measures used by Davini--Fathi--Iturriaga--Zavidovique.
In the linear case, it gives a short alternative proof of their
convergence and selection theorem. An example shows that linear and
cubic discounts may select critical solutions that do not differ by
an additive constant.
\end{abstract}


\maketitle
\enlargethispage{3pt}

\section{Introduction}

Let $M$ be a closed connected Riemannian manifold and let $H\in C(T^*M)$ be convex and superlinear in the momentum. We denote by $c_0$ the unique constant for which
\begin{equation}\label{critical}\tag{HJ$_0$}
H(x,Du)=c_0\qquad\text{in }M
\end{equation}
admits viscosity solutions. The discounted approximation was introduced in the work of
Lions, Papanicolaou and Varadhan \cite{LPV}. The convergence and
selection of solutions of
\[\lambda u_\lambda+H(x,Du_\lambda)=c_0\]
as $\lambda\to0^+$ were proved in \cite{DFIZ}.
The limit is the largest critical subsolution $w$ satisfying
\[\int_M w\,d\mu\le0\]
for every projected Mather measure $\mu$.
A variational approach using viscosity Mather measures for
second-order equations was developed in \cite{IMT}.

Nonlinear vanishing discount problems were studied in \cite{CCIZ,WYZ}. In particular, \cite{WYZ} treats smooth contact Hamiltonians satisfying $0<\partial_uG\le\Delta$. For the spatially degenerate equation \[\lambda\alpha(x)u_\lambda+H(x,Du_\lambda)=c_0,\] with $\alpha\ge0$ and $\alpha>0$ on the Aubry set, convergence was proved in \cite{Zav}. Convergence under an averaged nondegeneracy condition, $\int_{TM}\partial_uL_G(x,v,0)\,d\widetilde\mu<0$ for every Mather measure $\widetilde\mu$, was established in \cite{CFZZ}. Related convergence and divergence phenomena for non-monotone equations were studied in \cite{DNYZ}. For related results in the discrete setting, see \cite{DFIZdiscrete,NZdiscrete}.

The question of convergence for \[\lambda u_\lambda^3+H(x,Du_\lambda)=c_0\] was raised by Jinxin Xue. In this paper, we address this question and, more generally, nonlinear terms with real exponents. Although $u\mapsto u^3$ is strictly increasing, its derivative vanishes at $u=0$, so the preceding nondegeneracy conditions fail. This raises the question of whether convergence persists when
the first-order term vanishes, and how the leading nonlinear term
determines the selected critical solution. We prove convergence to the largest critical subsolution $w$ satisfying
\[\int_M w^3\,d\mu\le0\] for every projected Mather measure $\mu$,
in contrast to the condition $\int_M w\,d\mu\le0$
in the linear discount problem. We show in Section~\ref{example-section} that the two limits need not differ by a constant. For smooth nonlinearities, assumption (G3) below describes a common leading order in the Taylor expansion at the origin. More generally, it allows noninteger exponents without requiring differentiability in the last variable.

Different perturbations of the critical equation may lead to different selections of critical solutions. In \cite{NYZ}, sign-changing discount coefficients allow the selection of solutions associated with prescribed static classes. Other selection mechanisms arise by combining discount terms with potential perturbations \cite{ChenSelection,ChenZhu}, or by letting discount and viscosity vanish simultaneously \cite{MNT}. In contrast, the phenomenon studied here occurs for strictly increasing discount terms, without any additional perturbation: the order of vanishing of the discount term itself affects the selected solution.

More generally, let $G\in C(T^*M\times\R)$ and set $H(x,p):=G(x,p,0)$. We study
\begin{equation}\label{discount}\tag{E$_\lambda$}
G\bigl(x,Du_\lambda(x),\lambda u_\lambda(x)\bigr)=c_0\qquad\text{in }M,
\end{equation}
under the following assumptions:
\begin{itemize}
\item[(G1)] $p\mapsto G(x,p,r)$ is convex, and $G(x,p,r)/\|p\|_x\to+\infty$ as $\|p\|_x\to\infty$, uniformly for $x\in M$ and $|r|\le R$, for every $R>0$.
\item[(G2)] $r\mapsto G(x,p,r)$ is strictly increasing for every $(x,p)\in T^*M$.
\item[(G3)] For the Lagrangian \[L_G(x,v,r):=\sup_{p\in T_x^*M}\{p(v)-G(x,p,r)\},\] there exist a real number $q>0$ and $a\in C(TM)$, $a\ge0$, such that, for every compact $K\subset TM$,
\[
\bigl|L_G(x,v,r)-L_G(x,v,0)+a(x,v)f_q(r)\bigr|
\le |r|^q\omega_K(|r|),\qquad (x,v)\in K,\quad r\in\R,
\]
where $f_q(r):=\operatorname{sgn}(r)|r|^q$ (with $f_q(0)=0$), and $\omega_K:[0,\infty)\to[0,\infty)$ is nondecreasing and $\omega_K(\rho)\to0$ as $\rho\to0^+$.
\end{itemize}
The exponent $q$ and the function $a$ in (G3) are independent of $K$. We write $L_0:=L_G(\cdot,\cdot,0)$ and denote by $\MM$ the set of Mather measures for $L_0$, recalled in Section~\ref{prelim}. We further assume
\begin{equation}\label{positive}\tag{G4}
\int_{TM}a(x,v)\,d\widetilde\mu>0\qquad\text{for every }\widetilde\mu\in\MM.
\end{equation}
All solutions and subsolutions below are understood in the viscosity sense and are continuous.

\begin{theorem}\label{main}
Assume (G1)--(G4). For every $\lambda>0$, equation \eqref{discount} admits a unique solution $u_\lambda$. As $\lambda\to0^+$, $u_\lambda$ converges uniformly to the solution $u_0$ of \eqref{critical} given by
\[
u_0(x)=\sup_{w\in\mathcal S}w(x),\qquad
\mathcal S:=\left\{w:\ H(x,Dw)\le c_0,\quad
\int_{TM}a(x,v)f_q(w(x))\,d\widetilde\mu\le0\ \ \forall\widetilde\mu\in\MM\right\}.
\]
\end{theorem}

\begin{corollary}\label{power}
Let $H\in C(T^*M)$ be convex and superlinear in $p$, and let $b\in C(M)$, $b>0$. For every real number $q>0$, the unique solution of $\lambda b(x)f_q(u_\lambda)+H(x,Du_\lambda)=c_0$ converges uniformly as $\lambda\to0^+$ to
\[
\sup\left\{w(x):\ H(x,Dw)\le c_0,\quad
\int_M b(y)f_q(w(y))\,d\mu(y)\le0\ \ \forall\mu\in\mathfrak M\right\},
\]
where $\mathfrak M$ denotes the set of projected Mather measures of $H$.
\end{corollary}

\begin{proof}
Apply Theorem~\ref{main} to $G(x,p,r)=b(x)f_q(r)+H(x,p)$, with discount parameter $\lambda^{1/q}$.
\end{proof}

In particular, the result also covers $0<q<1$, where $f_q$ is not differentiable at the origin.

\begin{remark}\label{coercive-remark}
Taking $q=1$ and $b\equiv1$ in Corollary~\ref{power}, we recover
convergence and the selection formula of \cite{DFIZ} under the
superlinearity assumption. Our proof uses ordinary time averages along
truncated calibrated curves instead of the exponentially weighted
measures of \cite{DFIZ}. A related truncation argument along discounted calibrated curves
appears in the unpublished manuscript of Fathi and Iturriaga
\cite{FI}, where the truncation is based on level sets of the
difference between two critical solutions. For completeness, Appendix~\ref{app:DFIZ}
gives this argument under the original continuous, convex
and coercive assumptions of \cite{DFIZ}.

For separated Hamiltonians $G(x,p,r)=H(x,p)+V(x,r)$, with
$V(x,0)=0$, superlinearity can be replaced by uniform coercivity
in Corollary~\ref{power}.
Indeed, the uniform Lipschitz bounds allow us to modify $H$ outside
a fixed momentum ball and reduce to the superlinear case, as in
\cite{DFIZ} and \cite[Section~4.2]{Zav}. Since the discount term
is independent of the momentum, its expansion is unaffected by
this modification. For a general momentum-dependent $G$, preservation
of (G3) requires a separate justification.
\end{remark}

\begin{remark}
Suppose that $q$ is an odd positive integer and $L_G$ is $C^q$ in the last variable near $0$, locally uniformly on $TM$. If
$\partial_r^jL_G(x,v,0)=0$ for $1\le j<q$, then (G3) holds with
\[
a(x,v)=-\frac1{q!}\partial_r^q L_G(x,v,0),
\]
provided $a\ge0$. Condition (G4) must still be checked. The coefficient $a$ may vanish at some points, but (G3) requires a single exponent $q$ on every compact subset of $TM$.

The selection is determined by the leading term, rather than the full nonlinear dependence. For example, if
$G(x,p,r)=H(x,p)+b(x)(r+r^3)$ with $b>0$, then (G3) holds with $q=1$ and $a=b$, so the selection condition is $\int_M bw\,d\mu\le0$, not $\int_M b(w+w^3)\,d\mu\le0$.
\end{remark}

\begin{remark}
The theorem does not extend to very flat functions instead of $f_q$ as there are counterexamples with $f(r) = \operatorname{sgn}(r) \exp(-1/|r|)$. 
However, when the Aubry set is disconnected, the theorem can be
extended to functions $f_{q(x)}$ with different exponents $q(x)>0$
on different connected components of the Aubry set. We leave those developments to a future work in order to keep this one as elementary as possible. 
\end{remark}

Our proof of Theorem~\ref{main} compares critical subsolutions with
$u_\lambda$ along backward calibrated curves. The approaches of
\cite{WYZ,CFZZ} use exponentially weighted measures involving the
derivative of the Lagrangian in its last variable, extending the
construction of \cite{DFIZ}. Here, truncating the curves while the
difference remains positive allows us to use ordinary time averages
and the leading-order expansion directly, without differentiating
the Lagrangian in its last variable.

\par\medskip
\noindent\textbf{Organization of the paper.} Section~\ref{prelim} recalls the preliminary results.
Section~\ref{sec:main-proof} proves Theorem~\ref{main}.
Section~\ref{example-section} gives an example,
and Appendix~\ref{app:DFIZ} gives the truncation proof of the
classical linear discount theorem.

\section{Preliminaries}\label{prelim}

We recall the facts needed below; see \cite{Fathi}.

A function $w\in C(M)$ is a critical subsolution if and only if, for every absolutely continuous curve $\gamma:[s,t]\to M$,
\begin{equation}\label{domination}
w(\gamma(t))-w(\gamma(s))\le\int_s^t[L_0(\gamma(\tau),\dot\gamma(\tau))+c_0]\,d\tau.
\end{equation}
Critical subsolutions are equi-Lipschitz.

A Borel probability measure $\widetilde\mu$ on $TM$ is closed if \[\int_{TM}\|v\|_x\,d\widetilde\mu<\infty\quad \text{and} \quad \int_{TM}D\phi(x)[v]\,d\widetilde\mu=0\quad \text{for every } \phi\in C^1(M).\] If $\PP$ denotes the set of closed probability measures, then
\[
\min_{\widetilde\mu\in\PP}\int_{TM}L_0\,d\widetilde\mu=-c_0,
\qquad
\MM=\left\{\widetilde\mu\in\PP:\int_{TM}L_0\,d\widetilde\mu=-c_0\right\}.
\]
The set $\MM$ is nonempty and weakly$^*$ compact, and all its elements are supported in a common compact subset of $TM$. We write $\mu:=\pi_\#\widetilde\mu$, where $\pi:TM\to M$ is the projection.

A curve $\gamma:I\to M$ is calibrated by a critical subsolution $u$ if equality holds in \eqref{domination} with $w=u$ for all $s<t$ in $I$.

A function $u \in C(M)$ is a critical solution if and only if it is a critical subsolution and if for all $x\in M$, there is a calibrated curve with $I = (-\infty, 0]$ and $\gamma(0)=x$.

\section{Proof of Theorem~\ref{main}}\label{sec:main-proof}

Throughout this section, sequences of discount parameters satisfy $\lambda_j\to0^+$.

\subsection{Compactness}

By (G1)--(G2), the comparison principle holds for \eqref{discount}; see \cite[Theorem~3.2]{CCIZ}. Fix a critical solution $v_0$. The functions $v_0-\|v_0\|_\infty$ and $v_0+\|v_0\|_\infty$ are, respectively, a subsolution and a supersolution of \eqref{discount}. Perron's method (see \cite[Section~4]{CIL}) and comparison give a unique solution $u_\lambda$, with
\[
v_0-\|v_0\|_\infty\le u_\lambda\le v_0+\|v_0\|_\infty.
\]
Fix $C_0>0$ such that $\|u_\lambda\|_\infty\le C_0$. For $0<\lambda\le1$, (G2) gives $G(x,Du_\lambda,-C_0)\le c_0$ in the viscosity sense. By the coercivity of $p\mapsto G(x,p,-C_0)$, there exists $\ell>0$, independent of $\lambda$, such that $\Lip(u_\lambda)\le\ell$. Since $G(x,p,\lambda u_\lambda(x))\to H(x,p)$ locally uniformly on $T^*M$, every uniform limit of $u_\lambda$ as $\lambda\to0^+$ is a solution of \eqref{critical} by the stability of viscosity solutions.

Set $H_\lambda(x,p):=G(x,p,\lambda u_\lambda(x))$ and $L_\lambda(x,v):=L_G(x,v,\lambda u_\lambda(x))$. Since $H_\lambda(x,Du_\lambda)=c_0$, the critical value of $H_\lambda$ is $c_0$. For every $z\in M$, there exists a backward calibrated curve $\gamma_\lambda^z:(-\infty,0]\to M$, with $\gamma_\lambda^z(0)=z$, such that
\begin{equation}\label{calibration}
u_\lambda(\gamma_\lambda^z(t))-u_\lambda(\gamma_\lambda^z(s))
=\int_s^t[L_\lambda(\gamma_\lambda^z(\tau),\dot\gamma_\lambda^z(\tau))+c_0]\,d\tau,
\qquad s<t\le0.
\end{equation}
These curves are equi-Lipschitz. Indeed, $L_G$ is nonincreasing in $r$, and superlinearity gives \[L_\lambda(x,v)\ge L_G(x,v,C_0)\ge(\ell+1)\|v\|_x-C_1\] for some constant $C_1>0$. Differentiating \eqref{calibration} almost everywhere and using $\Lip(u_\lambda)\le\ell$ yields
\[
(\ell+1)\|\dot\gamma_\lambda^z\|+c_0-C_1
\le (u_\lambda\circ\gamma_\lambda^z)'\le\ell\|\dot\gamma_\lambda^z\|.
\]
Thus $\|\dot\gamma_\lambda^z\|\le C_1-c_0$ almost everywhere, uniformly in $z$ and $0<\lambda\le1$. Hence, all $(\gamma_\lambda^z,\dot\gamma_\lambda^z)$ lie almost everywhere in a compact disk bundle $K:=\{(x,v):\|v\|_x\le R\}$. We choose $R$ large enough that $K$ also contains the supports of all $\widetilde\mu\in\MM$.

\subsection{Mather measures}

By (G3), on $K$ we may write
\begin{equation}\label{expansion}
L_\lambda(x,v)=L_0(x,v)-\lambda^q a(x,v)f_q(u_\lambda(x))+R_\lambda(x,v),
\qquad |R_\lambda|\le\lambda^q\varepsilon_\lambda,
\end{equation}
where $\varepsilon_\lambda:=C_0^q\omega_K(\lambda C_0)\to0$.

\begin{lemma}\label{constraints}
If $u_{\lambda_j}\to u_*$ uniformly, then $u_*\in\mathcal S$.
\end{lemma}

\begin{proof}
Fix $\widetilde\mu\in\MM$. Since the critical value of $H_\lambda$ is $c_0$,
\[
\int_{TM}L_\lambda\,d\widetilde\mu
\ge\min_{\widetilde\nu\in\PP}\int_{TM}L_\lambda\,d\widetilde\nu
=-c_0=\int_{TM}L_0\,d\widetilde\mu.
\]
By \eqref{expansion}, \[\int_{TM}a(x,v)f_q(u_\lambda(x))\,d\widetilde\mu\le\varepsilon_\lambda.\] Passing to the limit gives the assertion.
\end{proof}

\begin{lemma}\label{curve-measure}
Let $\gamma_j:=\gamma_{\lambda_j}^{z_j}$, where $z_j\in M$, and let $T_j\to\infty$. Define the probability measures $\widetilde\nu_j$ by
\[
\int_{TM}f\,d\widetilde\nu_j:=\frac1{T_j}\int_{-T_j}^0f(\gamma_j(s),\dot\gamma_j(s))\,ds,\qquad f\in C_c(TM).
\]
Then every weak$^*$ accumulation point of $(\widetilde\nu_j)$ belongs to $\MM$.
\end{lemma}

\begin{proof}
All the measures are supported in $K$. Passing to a subsequence, assume $\widetilde\nu_j\rightharpoonup\widetilde\nu$. For every $\phi\in C^1(M)$,
\[
\int_{TM}D\phi(x)[v]\,d\widetilde\nu_j
=\frac{\phi(z_j)-\phi(\gamma_j(-T_j))}{T_j}\longrightarrow0.
\]
Since $K$ is compact, $\widetilde\nu$ is closed. By \eqref{calibration},
\[
\int_{TM}(L_0+c_0)\,d\widetilde\nu_j
=\frac{u_{\lambda_j}(z_j)-u_{\lambda_j}(\gamma_j(-T_j))}{T_j}
+\int_{TM}(L_0-L_{\lambda_j})\,d\widetilde\nu_j\longrightarrow0,
\]
because $\|u_{\lambda_j}\|_\infty\le C_0$ and $L_{\lambda_j}\to L_0$ uniformly on $K$. Thus $\int_{TM}(L_0+c_0)\,d\widetilde\nu=0$, so $\widetilde\nu\in\MM$.
\end{proof}

\subsection{Convergence}

Fix a sequence $\lambda_j\to0^+$ such that $u_{\lambda_j}\to u_*$ uniformly. Lemma~\ref{constraints} gives $u_*\in\mathcal S$. We prove $w\le u_*$ for every $w\in\mathcal S$.

Suppose otherwise. Choose $z\in M$ and $r>0$ such that $w(z)-u_*(z)>r$. Write $\gamma_j:=\gamma_{\lambda_j}^z$. By \eqref{domination}, \eqref{calibration}, and \eqref{expansion}, for every $T>0$,
\begin{equation}\label{comparison}
(w-u_{\lambda_j})(z)-(w-u_{\lambda_j})(\gamma_j(-T))\le\lambda_j^q\int_{-T}^0a(\gamma_j(s),\dot\gamma_j(s))f_q(u_{\lambda_j}(\gamma_j(s)))\,ds
+\lambda_j^q\varepsilon_{\lambda_j}T.
\end{equation}
For all sufficiently large $j$, $(w-u_{\lambda_j})(z)>r$. Define
\[
\tau_j:=\sup\{T>0:(w-u_{\lambda_j})(\gamma_j(s))>r\text{ for all }s\in(-T,0]\},
\qquad
T_j:=\begin{cases}\tau_j,&\tau_j<\infty,\\j/\lambda_j^q,&\tau_j=\infty.\end{cases}
\]
Then $(w-u_{\lambda_j})(\gamma_j(s))\ge r$ for $s\in[-T_j,0]$. If $\tau_j<\infty$, continuity gives $(w-u_{\lambda_j})(\gamma_j(-T_j))=r$, and \eqref{comparison} yields
\[
0<(w-u_{\lambda_j})(z)-r
\le\lambda_j^q\bigl(C_0^q\max_K a+\varepsilon_{\lambda_j}\bigr)T_j.
\]
Since the left-hand side converges to $w(z)-u_*(z)-r>0$, we have $T_j\to\infty$. This also holds when $\tau_j=\infty$.

Let $\widetilde\nu_j$ be the measures in Lemma~\ref{curve-measure} associated with $\gamma_j$ and $T_j$. Dividing \eqref{comparison} by $\lambda_j^qT_j$, we obtain
\[
\int_{TM}a(x,v)f_q(u_{\lambda_j}(x))\,d\widetilde\nu_j
\ge-\varepsilon_{\lambda_j}-\frac{2(\|w\|_\infty+C_0)}j.
\]
Indeed, the left-hand side of \eqref{comparison} is positive when $\tau_j<\infty$, and is bounded below by $-2(\|w\|_\infty+C_0)$ when $\tau_j=\infty$, in which case $\lambda_j^qT_j=j$.

Since $f_q$ is continuous and strictly increasing, for any
$C_w>\max\{\|w\|_\infty,C_0\}$ there exists $\delta>0$ such that
\[
f_q(t)-f_q(s)\ge\delta\qquad
\text{whenever }s,t\in[-C_w,C_w]\text{ and }t-s\ge r.
\]
Along $\gamma_j|_{[-T_j,0]}$ we have $w-u_{\lambda_j}\ge r$, and $a\ge0$, so
\[
\int_{TM}a(x,v)f_q(w(x))\,d\widetilde\nu_j
\ge\delta\int_{TM}a(x,v)\,d\widetilde\nu_j
-\varepsilon_{\lambda_j}-\frac{2(\|w\|_\infty+C_0)}j.
\]
By Lemma~\ref{curve-measure}, a subsequence converges to $\widetilde\nu\in\MM$. Passing to the limit and using \eqref{positive}, we obtain
\[
\int_{TM}a(x,v)f_q(w(x))\,d\widetilde\nu
\ge\delta\int_{TM}a(x,v)\,d\widetilde\nu>0,
\]
contrary to $w\in\mathcal S$.

We have proved $w\le u_*$ for every $w\in\mathcal S$. Together with $u_*\in\mathcal S$, this gives $u_*(x)=\sup_{w\in\mathcal S}w(x)$. Every convergent subsequence therefore has the same limit, and compactness yields uniform convergence of the whole family. This proves Theorem~\ref{main}.

\begin{remark}
The stopping times in this proof cannot be chosen freely, since
$w-u_{\lambda_j}\ge r$ must hold along each truncated curve.
Consequently, information obtained from infinite-time averages in
the non-monotone problem of \cite{NYZ} cannot be applied directly
to these segments. In particular, the argument does not remove
the hypothesis on the Mather quotient in \cite{NYZ}.
\end{remark}

\section{An example}\label{example-section}

We show that the linear and cubic discounts may select critical solutions that do not differ by a constant.
Let $M=\mathbb T^2=\mathbb R^2/\mathbb Z^2$, with coordinates $(x,y)$, and set
\[
f(x)=\cos(2\pi x)+\frac12\cos(4\pi x),\qquad
b(y)=\frac{1-\cos(2\pi y)}2,\qquad
\psi(x,y)=b(y)f(x).
\]
Consider the Tonelli Hamiltonian
\[
H(x,y,p_x,p_y)=\frac12\bigl(p_x+1+\partial_x\psi(x,y)\bigr)^2
+\frac12\bigl(p_y+\partial_y\psi(x,y)\bigr)^2
+\cos^2(2\pi y).
\]
Its Lagrangian is
\[
L(x,y,v_x,v_y)=\frac12(v_x^2+v_y^2)-(1+\partial_x\psi)v_x
-\partial_y\psi\,v_y-\cos^2(2\pi y).
\]
In particular,
\begin{equation}\label{example-action}
L+\frac32=\frac12(v_x-1)^2+\frac12v_y^2
+\sin^2(2\pi y)-D\psi(x,y)\cdot(v_x,v_y).
\end{equation}
The last term has zero integral against every closed measure. Hence
$\int(L+3/2)\,d\widetilde\mu\ge0$ for every closed probability
measure $\widetilde\mu$. Equality holds for
\[
\widetilde\mu_0=dx\otimes\delta_{y=0}\otimes\delta_{(1,0)},
\qquad
\widetilde\mu_1=dx\otimes\delta_{y=1/2}\otimes\delta_{(1,0)},
\]
where $dx$ is the normalized Lebesgue measure on $\mathbb T^1$.
Hence, $c_0=3/2$. Moreover, equality in \eqref{example-action} forces
$v_x=1$, $v_y=0$ and $\sin(2\pi y)=0$ almost everywhere.
The closedness condition then forces the $x$-marginal on each of the
two circles $y=0$ and $y=1/2$ to be Lebesgue measure.
Consequently, the Mather measures are precisely the convex combinations
of $\widetilde\mu_0$ and $\widetilde\mu_1$.

Let $w$ be any critical subsolution. Along either periodic curve
$\gamma_i(t)=(x+t,i/2)$, $i=0,1$, the right-hand side of
\eqref{domination} is
\[
\int_s^t(L(\gamma_i,\dot\gamma_i)+c_0)\,d\tau
=-\psi(\gamma_i(t))+\psi(\gamma_i(s)).
\]
Thus $w+\psi$ is nonincreasing along each of these periodic curves.
By periodicity it must be constant on each circle. Since
$\psi(x,0)=0$ and $\psi(x,1/2)=f(x)$, there exist
$\alpha_0,\alpha_1\in\mathbb R$ such that
\begin{equation}\label{example-values}
w(x,0)=\alpha_0,\qquad
w(x,1/2)=\alpha_1-f(x).
\end{equation}

Consider first the linear discount equation
$\lambda u_\lambda+H(x,Du_\lambda)=c_0$, and denote its limit
by $u_0^{(1)}$. By Corollary~\ref{power}, the selection constraints
for $w$ are $\int w\,d\mu_i\le0$, $i=0,1$, where
$\mu_i$ is the projection of $\widetilde\mu_i$.
Since $\int_0^1 f(x)\,dx=0$, \eqref{example-values} gives
\[
\alpha_0\le0,\qquad\alpha_1\le0.
\]
On the other hand, $w_1=-\psi$ is a critical subsolution, since
\[
H(x,y,Dw_1)=\frac12+\cos^2(2\pi y)\le c_0,
\]
and its integrals against both $\mu_i$ vanish. By the maximality
in Corollary~\ref{power}, we obtain
\begin{equation}\label{example-linear}
u_0^{(1)}(x,0)=0,\qquad
u_0^{(1)}(x,1/2)=-f(x).
\end{equation}

Next, consider the cubic discount equation
$\lambda u_\lambda^3+H(x,Du_\lambda)=c_0$, and denote its limit
by $u_0^{(3)}$. The corresponding selection constraints are
$\int w^3\,d\mu_i\le0$ for $i=0,1$. A direct calculation gives
\[
\int_0^1 f(x)^2\,dx=\frac58,
\qquad
\int_0^1 f(x)^3\,dx=\frac38.
\]
In view of \eqref{example-values}, the constraints become
\[
\alpha_0\le0,\qquad
\alpha_1^3+\frac{15}{8}\alpha_1-\frac38\le0.
\]
Let $\kappa$ be the unique root of
\begin{equation}\label{example-kappa}
\kappa^3+\frac{15}{8}\kappa-\frac38=0.
\end{equation}
The polynomial is strictly increasing, and its values at $0$ and
$1/4$ have opposite signs. In particular, $0<\kappa<1/4$, and
$\alpha_1\le\kappa$.

To see that both bounds can be attained, put
$w_3=-\psi+\kappa b(y)$. Since $b'(y)=\pi\sin(2\pi y)$,
\[
H(x,y,Dw_3)
=\frac12+\frac{\kappa^2\pi^2}{2}\sin^2(2\pi y)
+\cos^2(2\pi y) =\frac32-
\left(1-\frac{\kappa^2\pi^2}{2}\right)
\sin^2(2\pi y) \le\frac32=c_0,
\]
where the last inequality follows from
$\kappa^2\pi^2/2<\pi^2/32<1$.
Thus $w_3$ is a critical subsolution. It has values
$w_3(x,0)=0$ and $w_3(x,1/2)=\kappa-f(x)$, so its cubic integrals
against both Mather measures vanish, by \eqref{example-kappa}.
Theorem~\ref{main} (or Corollary~\ref{power}) and maximality now give
\begin{equation}\label{example-cubic}
u_0^{(3)}(x,0)=0,\qquad
u_0^{(3)}(x,1/2)=\kappa-f(x).
\end{equation}
Comparing \eqref{example-linear} and \eqref{example-cubic}, we conclude that
\[
(u_0^{(3)}-u_0^{(1)})(x,0)=0,
\qquad
(u_0^{(3)}-u_0^{(1)})(x,1/2)=\kappa>0.
\]
In particular, the two selected critical solutions do not differ by
an additive constant.

\appendix
\section{The linear discount problem}\label{app:DFIZ}

We record the truncation argument for the classical equation under
the coercivity assumption of \cite{DFIZ}. This proof does not use
the exponentially weighted measures introduced there.

\begin{proposition}
Let $H\in C(T^*M)$ be convex and coercive in the momentum, uniformly
in $x$. Let $u_\lambda$ solve
\[
 \lambda u_\lambda+H(x,Du_\lambda)=c_0.
\]
Then $u_\lambda$ converges uniformly as $\lambda\to0^+$ to the largest
critical subsolution $w$ such that
\[
 \int_M w\,\dd\mu\le0
 \quad\text{for every projected Mather measure }\mu.
\]
\end{proposition}

\begin{proof}
Comparison with translates of a critical solution gives a uniform
bound for $u_\lambda$. Coercivity then gives a common Lipschitz bound
$\Lip(u_\lambda)\le\ell$ for $0<\lambda\le1$.
By the standard modification of $H$ outside a fixed momentum ball,
as in \cite{DFIZ} and Remark~\ref{coercive-remark}, we may assume
that $H$ is superlinear. Let $L$ be its Lagrangian. The backward
calibrated curves $\gamma_\lambda$ for $u_\lambda$ have uniformly
bounded velocities, by the compactness argument in
Section~\ref{sec:main-proof}.

Suppose $\lambda_j\to0$ and $u_{\lambda_j}\to u_*$ uniformly.
The limit $u_*$ is a critical solution. For any Mather measure
$\widetilde\mu$ with projection $\mu$, the Hamiltonian
$H+\lambda_j u_{\lambda_j}(x)$ has critical value $c_0$.
Its Lagrangian is $L-\lambda_j u_{\lambda_j}(x)$, and hence
\[
 -c_0\le\int_{TM}(L-\lambda_j u_{\lambda_j})\,\dd\widetilde\mu
 =-c_0-\lambda_j\int_Mu_{\lambda_j}\,\dd\mu.
\]
Thus $\int_Mu_*\,\dd\mu\le0$.

Let $w$ be any critical subsolution satisfying these inequalities.
Assume that $w(z)>u_*(z)$ for some $z\in M$. Choose $r>0$ such that
$w(z)-u_*(z)>r$, and let $\gamma_j:(-\infty,0]\to M$ be a
backward calibrated curve for $u_{\lambda_j}$ ending at $z$.
Domination for $w$ and calibration for $u_{\lambda_j}$ give
\[
 (w-u_{\lambda_j})(z)
 -(w-u_{\lambda_j})(\gamma_j(-T))
 \le\lambda_j\int_{-T}^0u_{\lambda_j}(\gamma_j(s))\,\dd s.
\]
Let $\tau_j$ be the first backward time at which the difference
$w-u_{\lambda_j}$ reaches $r$, and put $T_j=\tau_j$ if this time is
finite and $T_j=j/\lambda_j$ otherwise. On $[-T_j,0]$ we have
$w-u_{\lambda_j}\ge r$. In the finite case, the positive difference
$(w-u_{\lambda_j})(z)-r$ and the uniform bound on $u_{\lambda_j}$
imply $\lambda_j T_j\ge c>0$. In the infinite case,
$\lambda_jT_j=j$. Thus $T_j\to\infty$.

Consider the ordinary time averages
\[
 \int_{TM}\zeta\,\dd\widetilde\nu_j
 =\frac1{T_j}\int_{-T_j}^0
 \zeta(\gamma_j(s),\dot\gamma_j(s))\,\dd s,
 \qquad \zeta\in C_c(TM).
\]
Their supports lie in a fixed compact set. For any $\psi\in C^1(M)$,
\[
 \int_{TM}D\psi(x)[v]\,\dd\widetilde\nu_j
 =\frac{\psi(z)-\psi(\gamma_j(-T_j))}{T_j}\longrightarrow0.
\]
Moreover, calibration yields
\[
 \int_{TM}(L+c_0)\,\dd\widetilde\nu_j
 =\frac{u_{\lambda_j}(z)-u_{\lambda_j}(\gamma_j(-T_j))}{T_j}
  +\lambda_j\int_{TM}u_{\lambda_j}\,\dd\widetilde\nu_j
 \longrightarrow0.
\]
Since $L$ is lower semicontinuous and bounded below on this compact
set, every weak limit $\widetilde\nu$ is closed and satisfies
$\int(L+c_0)\,\dd\widetilde\nu\le0$. The minimal-action
characterization of $c_0$ gives the opposite inequality for a closed
probability measure. Hence $\widetilde\nu$ is a Mather measure.

Divide the gap comparison above by $\lambda_jT_j$. In the finite
case its left-hand side is positive; in the infinite case it is
bounded below independently of $j$, while $\lambda_jT_j=j$.
Consequently,
\[
 \int_{TM}u_{\lambda_j}(x)\,\dd\widetilde\nu_j\ge-o(1).
\]
Since $w-u_{\lambda_j}\ge r$ along the selected segments,
\[
 \int_{TM}w(x)\,\dd\widetilde\nu_j\ge r-o(1).
\]
Passing to a Mather limit gives $\int_Mw\,\dd\mu\ge r>0$,
contrary to the condition on $w$. Thus $w\le u_*$.
As $u_*$ itself satisfies the integral constraints, it is the
largest admissible critical subsolution. Every convergent subsequence
has the same limit, which proves uniform convergence.
\end{proof}

\section*{Acknowledgements}
Panrui Ni is supported by JSPS KAKENHI Grant No.~26KF0103 and the National Natural Science Foundation of China (Grant No.~12571197). Maxime Zavidovique is supported by ANR CoSyDy (ANR-CE40-0014). M.Z. gratefully acknowledges support from the Institute for Advanced Study, of which he was a member when this work was finished. He is grateful to the anonymous donor funding his stay there.

\section*{Declarations}
\noindent\textbf{Conflict of interest.} The authors declare that they have no conflict of interest.

\smallskip
\noindent\textbf{Data availability.} Data sharing is not applicable to this article as no datasets were generated or analysed during the current study.

\smallskip
\noindent\textbf{AI.} All of the present mathematics is purely human made with no use of AI.

\end{document}